\documentstyle[
epsfig,verbatim
]{amsart}

\newcommand \lbb[1] {\label{#1}}

\newcommand {\1}    {{\bf 1}}

\newcommand     {\CD}   {{\C[\partial]}}
\newcommand     {\D}    {{\partial}}

\newcommand     {\C}    {{\Bbb C}}
\newcommand     {\N}    {{\Bbb N}}
\newcommand     {\Z}    {{\Bbb Z}}

\newcommand     {\End}  {\mbox{End}}

\newcommand     {\conf} {{\mbox{\scriptsize {\it Lie}}}}

\newcommand     {\id}   {\mbox{id}}
\newcommand     {\open} {``}
\newcommand     {\close}{''}

\newtheorem{thm}{Theorem}[section]

\newtheorem{lemma}{Lemma}[section]
\newtheorem{prop}{Proposition}[section]
\newtheorem{cor}{Corollary}[section]

\theoremstyle{definition}

\newtheorem{defn}{Definition}[section]

\newtheorem{rem}{Remark}[section]

\begin{document}

\title[Simplicity of Vertex and Lie Conformal Algebras]
{A remark on simplicity of vertex algebras and Lie conformal
algebras}

\author[A.~D'Andrea]{Alessandro D'Andrea}
\thanks{The author was partially supported by
PRIN ``Spazi di Moduli e Teoria di Lie'' fundings from MIUR and
project MRTN-CT 2003-505078 ``LieGrits'' of the European Union}
\address{Dipartimento di Matematica, Universit\`a degli Studi di
Roma ``La Sapienza''\\ P.le Aldo Moro, 5 -- 00185 Rome, Italy}
\email{dandrea@@mat.uniroma1.it}

%\dedicatory{Ad Olivia, che mi ha insegnato a salutare il Sole}

\begin{abstract}
I give a short proof of the following algebraic statement: if $V$ is
a simple vertex algebra, then the underlying Lie conformal algebra
is either abelian, or it is an irreducible central extension of a
simple Lie conformal algebra. This provides many examples of
non-finite simple Lie conformal algebras, and should prove useful
for classification purposes.
\end{abstract}

%\date{\today}

\maketitle

%\tableofcontents

\section{Introduction}

The notion of {\em vertex algebra} was introduced by Borcherds in
\cite{Bo} to axiomatize algebraic properties of the Operator Product
Expansion (= OPE) of quantum fields in a (chiral) Conformal Field
Theory in dimension two. Vertex algebras were defined as vector
spaces endowed with infinitely many bilinear operations satisfying
involved axioms that are now known as Borcherds identities.

The construction of non-trivial examples of vertex algebras is a
complicated matter, because interesting known examples are very
large objects -- typically graded vector spaces of superpolynomial
growth, called {\em Vertex Operator Algebras} -- and
finite-dimensional instances degenerate into differential
commutative algebra structures \cite{Bo}. Examples of physical
interest are usually described by giving generating quantum fields,
after prescribing the singular part of their OPE. This idea can be
made precise by axiomatizing the singular part of the OPE into a
{\em (Lie) conformal algebra} structure, introduced by Kac in
\cite{K}. Lie conformal algebras \cite{DK}, and their
generalizations \cite{DsK}, only determine commutation properties of
quantum fields, and the whole vertex algebra can then be recovered
by taking a suitable quotient of a certain universal envelope
\cite{K,Li,P,R} of the Lie conformal algebra.

The Lie conformal algebra theory has proved simpler than the vertex
algebra one. On the one hand, it is easy to construct small
non-trivial examples; on the other hand, Lie conformal algebras
possess a close resemblance to Lie algebras -- they are indeed Lie
algebras in a suitable pseudo-tensor category, see \cite{BKV,BDK} --
and can be treated by means of similar techniques.

It is clear that every vertex algebra is also a Lie conformal
algebra: the resulting forgetful functor is adjoint to the
above-mentioned universal enveloping vertex algebra functor. Both
vertex algebras and Lie conformal algebras have corresponding
notions of ideal and simplicity; however, it is easier for a
subspace to be an ideal with respect to the Lie conformal algebra
structure than with respect to the vertex one. The main result of
the present paper is a short and elementary proof of the quite
surprising fact that the Lie conformal algebra structure underlying
a simple vertex algebra is as simple as it can be: its only ideals
are central, and the whole Lie conformal algebra is a central
extension of a simple structure. Indeed the one-dimensional vector
space spanned by the vacuum element is always a central Lie
conformal ideal.

The main tool employed in the paper is identity \eqref{genwickeq}
whose constant (in $z$) part generalizes a formula devised by Wick
\cite{Wick} to compute the singular OPE of normally ordered products
of fields in a free theory, which was independently mentioned in
\cite{BK} and \cite{preprint}, and whose algebraic consequences
range beyond the present result.

\section{Vertex algebras}
Let $V$ be a complex vector space. A {\em field} on $V$ is a formal
distribution $\phi \in (\End V)[[z, z^{-1}]]$ with the property that
$\phi(v) \in V((z))$ for every $v \in V$. In other words, if
$$\phi(z) = \sum_{i \in \Z} \phi_i \, z^{-i-1},$$ then $\phi_n(v) = 0 $
for sufficiently large $n$.

\begin{defn}[\cite{K}]
A vertex (super)algebra is a complex vector super-space $V = V_0
\oplus V_1$ endowed with a linear parity preserving {\em state-field
correspondence} \mbox{$Y:V \to (\End V)[[z, z^{-1}]]$}, a {\em
vacuum element} $\1\in V_0$ and an even operator \mbox{$T \in \End
V$} satisfying the following properties:
\begin{itemize}
\item {\bf Field axiom}: $Y(v,z)$ is a field for all $v\in V$.
\item {\bf Locality}: For every choice of $a\in V_{p(a)}$ and $b \in V_{p(b)}$ one has
$$(z-w)^N \left(Y(a,z) Y(b,w) - (-1)^{p(a)p(b)} Y(b,w) Y(a,z)\right) = 0$$ for sufficiently large $N$.
\item {\bf Vacuum axiom}: The vacuum element $\1$ is such that
$$Y(\1,z) = \id_V,\qquad Y(a,z)\1 \equiv a \mod
zV[[z]],$$ for all $a \in V$.
\item {\bf Translation invariance}: $T$ satisfies $$[T, Y(a,z)]
= Y(Ta, z) = \frac{d}{dz}Y(a,z),$$ for all $a\in V$.
\end{itemize}
\end{defn}
Notice that the vector space $V$ carries a natural $\C[T]$-module
structure. Fields $Y(a,z)$ are called {\em vertex operators}, or
{\em quantum fields}. Vertex algebra axioms have strong algebraic
consequences, among which we recall the following:

\begin{itemize}
\item {\bf Skew-commutativity}: \qquad $Y(a,z)b = e^{zT} Y(b,-z)a$.
%\item {\bf Associativity}: \qquad $Y(a,z)Y(b,w) = Y(Y(a,z-w)b,w)$ in
%the domain $|z|>|w|$.
\end{itemize}
%Equality \open in the domain $|z|>|w|$\close\, means that before
%comparing the two sides, negative powers of $z-w$ in the expression
%$Y(Y(a,z-w)b,w)$ should be expanded as power series in $z$ and $w$
%assuming $|w|/|z|
%>1$. For instance
%$$\frac{1}{z-w} = \frac{1}{z}\cdot\frac{1}{1 - w/z} = \frac{1}{z}
%\cdot \sum_{k=0}^\infty \left(\frac{w}{z}\right)^k.$$
%\begin{rem}\lbb{axioms2}
%Commutativity and associativity of quantum fields, together with the
%vacuum axiom and translation invariance, provide an equivalent set
%of axioms for a vertex algebra \cite{Li}.
%\end{rem}
Coefficients of quantum fields
$$Y(a, z) = \sum_{j\in \Z} a_{(j)} z^{-j-1}$$
in a vertex algebra span a Lie algebra under the commutator Lie
bracket, and more explicitly satisfy -- see \cite[Theorem
2.3(iv)]{K},
\begin{equation}\lbb{bracket}
[a_{(m)}, b_{(n)}] = \sum_{j\in \N} {m\choose j}
(a_{(j)}b)_{(m+n-j)},
\end{equation}
for all $a,b\in V, m,n \in \Z$.

If $A$ and $B$ are subspaces of $V$, then we may define $A \cdot B$
as the $\C$-linear span of all products $a_{(j)} b$, where $a\in A,
b \in B, j \in \Z$. It follows that if $A$ and $B$ are
$\C[T]$-submodules of $V$, then $A\cdot B$ is also a
$\C[T]$-submodule of $V$, as by translation invariance $T$ is a
derivation of all $j$-products. Notice that in this case $A \cdot B
= B \cdot A$ by skew-commutativity, and that $A \subset A \cdot V$
by the vacuum axiom.

An {\em ideal} of $V$ is a $\C[T]$-submodule $I\subset V$ such that
$V \cdot I \subset I$. We will say that a vertex algebra is {\em
simple} if its only ideals are trivial.

\begin{rem}\lbb{idealcheck}
In order for a subspace $A \subset V$ to be an ideal, it suffices to
check that $A \cdot V \subset A$: then $A$ is indeed a
$\C[T]$-submodule of $V$, as $a_{(-2)} \1 = T a$; moreover,
skew-commutativity gives $V \cdot A = A \cdot V \subset A$.
\end{rem}

\section{Conformal algebras}

\begin{defn}[\cite{DK}] A {\em Lie conformal (super)algebra} is a
$\Z/2\Z$-graded $\C[\partial]$-module $R = R_0 \oplus R_1$ endowed
with a parity preserving $\C$-bilinear product $(a,b)\mapsto
[a\,_\lambda\, b] \in R[\lambda]$ satisfying the following axioms:
\begin{enumerate}
\item[{\bf (C1)}] $[a \,_\lambda\, b] \in R[\lambda],$

\item[{\bf (C2)}] $[\partial a \,_\lambda\, b] = -\lambda [a \,_\lambda\,
b],\,\,\, [a \,_\lambda\, \partial b] = (\partial +\lambda)
[a\,_\lambda\,b],$

\item[{\bf (C3)}] $[a \,_\lambda\, b] = - (-1)^{p(a)p(b)} [b \,_{-\partial -\lambda}\, a],$

\item[{\bf (C4)}] $[a \,_\lambda\, [b \,_\mu\, c]] - (-1)^{p(a)p(b)} [b \,_\mu\,
[a\,_\lambda\,c]] = [[a \,_\lambda\, b] \,_{\lambda + \mu}\, c],$
\end{enumerate}
\noindent for every choices of homogeneous elements $a, b, c \in V$,
$p(r)\in \Z/2\Z$ denoting the parity of the homogeneous element $r$.
\end{defn}

Any vertex (super)algebra $V$ can be given a $\CD$-module structure
by setting $\D = T$. Then defining
$$ [a\,_\lambda\, b] = \sum_{n\in \N} \frac{\lambda^n}{n!}
a_{(n)}b$$ endows $V$ with a Lie conformal (super)algebra structure.
Indeed, (C1) follows from the field axiom, (C2) from translation
invariance, (C3) from skew-commutativity, and (C4) from
\eqref{bracket}. For the sake of simplicity, in all that follows the
super- prefix will not be explicitly mentioned, but tacitly
understood.

If $A$ and $B$ are subspaces of a Lie conformal algebra $R$, then we
may define $[A ,B]$ as the $\C$-linear span of all
$\lambda$-coefficients in the products $[a\,_\lambda b]$, where
$a\in A, b \in B$. It follows from axiom (C2) that if $A$ and $B$
are $\CD$-submodules of $R$, then $[A, B]$ is also a $\CD$-submodule
of $R$. Notice that in this case $[A, B] = [B, A]$ by axiom (C3). A
Lie conformal algebra $R$ is {\em solvable} if, after defining
$$R^0 = R,\qquad R^{n+1} = [R^n, R^n], n\geq 0,$$
we find that $R^N = 0$ for sufficiently large $N$. $R$ is solvable
iff it contains a solvable ideal $S$ such that $R/S$ is again
solvable. Solvability of a nonzero Lie conformal algebra $R$
trivially fails if $R$ equals its derived subalgebra $R' = [R, R]$.
An {\em ideal} of a Lie conformal algebra $R$ is a $\CD$-submodule
$I\subset R$ such that $[R, I] \subset I$. If $I, J$ are ideals of
$R$, then $[I, J]$ is an ideal as well. An ideal $I$ is said to be
{\em central} if $[R, I] = 0$, i.e., if it is contained in the {\em
centre} $Z=\{r \in R| [r\,_\lambda s] = 0 \mbox{ for all } s \in
R\}$ of $R$. A Lie conformal algebra $R$ is {\em simple} if its only
ideals are trivial, and $R$ is not {\em abelian}, i.e., $[R, R] \neq
0$.

Notice that, when $V$ is a vertex algebra, we should distinguish
between ideals of the vertex algebra structure and ideals of the
underlying Lie conformal algebra. Indeed, ideals of the vertex
algebra are also ideals of the Lie conformal algebra, but the
converse is generally false, as it can be seen by observing that
$\C\1$ is always a central ideal of the Lie conformal algebra
structure, but it is never an ideal of the vertex algebra.

In order to avoid confusion, we will denote by $V^\conf$ the Lie
conformal algebra structure underlying a vertex algebra $V$.

\section{A Poisson-like generalization of the Wick formula}

The following formula -- which is similar to (3.3.7) and (3.3.12) in
\cite{K} -- relating the vertex and Lie conformal algebra structures
is the key tool in the present paper.
\begin{prop}\lbb{genwick}
If $a,b,c$ are elements of the vertex algebra $V$, then:
\begin{equation}\lbb{genwickeq}
[a\,_\lambda \,Y(b, z)\,c] = e^{\lambda z} Y([a\,_\lambda\, b],
z)\,c + Y(b, z) \,[a\,_\lambda \,c].
\end{equation}
\end{prop}
\begin{pf}
Multiply both sides of \eqref{bracket} by $\lambda^m z^{-n-1}/m!$,
then add up over all $m\in \N, n\in \Z$. Applying both sides to
$c\in V$ proves the statement.
\end{pf}

\begin{lemma}\lbb{key}
Let $U \subset V$ be vector spaces, and $p(\lambda), q(\lambda)$ be
elements of $V((z))[\lambda]$. If all coefficients of
$$e^{\lambda z} p(\lambda) + q(\lambda)$$
lie in $U((z))$, then the same is true for the coefficients of
$p(\lambda)$.
\end{lemma}
\begin{pf}
If $m$ and $n$ are the degrees of $p$ and $q$ as polynomials in
$\lambda$, we write
$$p(\lambda) = \sum_{i=0}^m p_i(z) \lambda^i,\qquad q(\lambda) =
\sum_{j=0}^n q_i(z) \lambda^j.$$ The expression $e^{\lambda z}
p(\lambda) + q(\lambda)$ is a power series in $\lambda$, and the
coefficient multiplying $\lambda^N$ is independent of $q(\lambda)$
if $N>n$. If also $N \geq m$, it equals
$$\sum_{i=0}^m \frac{z^{N-i}}{(N-i)!} p_i(z)= z^N \sum_{i=0}^m \frac{1}{(N-i)!}
\cdot \frac{p_i(z)}{z^i}.$$ If all $\lambda$-coefficients of
$e^{\lambda z} p(\lambda) + q(\lambda)$ lie in $U((z))$, then
$$\sum_{i=0}^m \frac{1}{(N-i)!} \cdot \frac{p_i(z)}{z^i} \in U((z))$$
for all sufficiently large $N$. However, the $(m+1) \times (m+1)$
matrix whose $(i, j)$-entry is $1/(N+j-i)!$ is
non-singular\footnote{Its determinant can be computed by induction,
and is easily showed to be equal to \mbox{$m! (m-1)! (m-2)! ... 3!
2! 1!/((N+m)! (N+m-1)! ... N!)$}.}, hence
$$\frac{p_i(z)}{z^i} \in U((z)),$$ thus showing that
$p(\lambda) \in U((z))[\lambda]$.
\end{pf}

\section{A simplicity argument}

\begin{thm}
Let $V$ be a vertex algebra, and $I \subset V$ a subspace. Then $[I,
V]$ is an ideal of $V$.
\end{thm}
\begin{pf}
Choose $a \in I, b,c \in V$. The linear span of all coefficients of
$[a_\lambda Y(b, z)c]$, when $a\in I, b, c\in V$ equals $[I, V \cdot
V] = [I, V]$. By Lemma \ref{key} applied to \eqref{genwickeq}, all
coefficients of $Y([a_\lambda b], z)c, a\in I, b,c\in V$ lie in $[I,
V]$, thus $[I, V]\cdot V \subset [I, V]$. Then Remark
\ref{idealcheck} ensures that $[I, V]$ is an ideal of $V$.
%This shows that $[I, V]$ is a
%$\C[T]$-submodule of $V$, and is moreover a right, hence two-sided,
%ideal of $V$.
\end{pf}

\begin{cor}\lbb{simple}
Let $V$ be a simple vertex algebra. Then either $V^\conf$ is
abelian, or it is an irreducible central extension of a simple Lie
conformal algebra.
\end{cor}
\begin{pf}
Let $I$ be a proper ideal of $V^\conf$. Then $[I, V] \subset I$ is a
proper ideal of $V$, forcing $[I, V] = 0$ by simplicity. Thus all
proper ideals of $V^\conf$ lie in the centre $Z$ of $V^\conf$, hence
either $V^\conf = Z$, or $V^\conf/Z$ has no non-trivial ideal; in
the former case $V^\conf$ is abelian.

In the latter, $[V, V]$ is a nonzero ideal of $V$, hence $V=[V, V]$.
Then $V^\conf$ is not solvable, as it equals its derived subalgebra,
so $V^\conf/Z$ cannot be abelian, and is therefore simple. As
$V^\conf$ equals its derived algebra, it is an irreducible central
extension.
\end{pf}

\begin{rem}
Corollary \ref{simple} can be used, along with the classification
\cite{DK} of finite simple (purely even) Lie conformal algebras and
a known characterization of irreducible central extensions of the
Virasoro Lie conformal algebra, to show that all simple (purely
even) vertex algebra structures on a finitely generated $\CD$-module
are {\em abelian}, i.e. have a trivial underlying Lie conformal
algebra structure. From this it follows that if $V$ is a finite
vertex algebra, then $V^\conf$ is solvable.

A more detailed investigation of such finite vertex algebras can
show that $V^\conf$ is indeed nilpotent, as soon as $V$ contains no
element $a$ such that $Y(a,z)a = 0.$\footnote{This is analogous to
demanding that a commutative algebra possess no nilpotent elements.}
Such claims are proved in a separate paper \cite{nilpotence}.
\end{rem}

\begin{rem}
Let $V$ be a vertex algebra, and assume that whenever a subspace
$U\subset V$ is invariant under the action of coefficients of all
quantum fields, then $U$ is a $\CD$-module, and therefore an ideal.
This happens, for instance, if $\D = T$ is a coefficient of some
quantum field, e.g., in a (conformal) Vertex Operator Algebra, where
$T=L_{-1}$.

By a Schur Lemma argument, one may then show that if the vertex
algebra $V$ is simple and countable-dimensional -- as it is always
the case when $V$ is a $\Z$-graded vector space with finite
dimensional homogeneous components -- then the only central elements
in the underlying Lie conformal algebra are scalar multiples of the
vacuum element. Then $V^\conf$ is an irreducible central extension
of a simple Lie conformal algebra by the one-dimensional ideal
$\C\1$.
\end{rem}

\begin{rem}
If $R$ is an irreducible central extension of a simple Lie conformal
algebra by a one-dimensional centre $\C\1$, and a grading is given
on $R$ which is compatible with its Lie conformal algebra structure,
then there exists at most one simple vertex algebra structure
compatible with the same grading, in which $\1$ is the vacuum
element. Indeed, the universal enveloping vertex algebra of $R$ has
a unique maximal graded ideal, which must intersect $R$ trivially
because of the vacuum axiom, and of the simplicity of $R/\C\1$.

This provides a strategy for finding simple vertex algebras, by
first looking for simple Lie conformal algebra structures and their
possible central extensions $R$, and then checking whether the
unique simple quotient of the universal enveloping vertex algebra is
$R$ or a larger space. This strategy might become effective for
families of Lie conformal algebras for which a classification of
simple objects is likely to be achieved, e.g., under a polynomial
growth or a finite Gelfand-Kirillov dimension assumption \cite{Xu},
\cite{Z1, Z2}. Recall that no Vertex Operator Algebra of physical
interest is of this kind, as the presence of a Virasoro field forces
superpolynomial growth.
\end{rem}

\end{document}